\documentclass{amsart}

\usepackage{amsmath} 
\usepackage{amssymb}

\newcommand{\forces}{\Vdash} 

\newcommand{\incomp}{\bot}
\newcommand{\compatible}{\mbox{$\not\!\bot$}}
 
\newcommand{\lbv}{[\![} 
\newcommand{\rbv}{]\!]}
 
\newcommand{\lesdot}{\mathrel{\mathord{<}\!\!\raise 
0.8 pt\hbox{$\scriptstyle\circ$}}} 
\newcommand{\comp}{\circ} 
  
\newcommand{\con}{{\mathfrak c}}

\newcommand{\can}{2^{\textstyle \omega}} 
\newcommand{\fs}{2^{\textstyle <\!\omega}} 
\newcommand{\baire}{\omega^{\textstyle \omega}} 
\newcommand{\iso}{[\omega]^{\textstyle \omega}} 
\newcommand{\fsuo}{[\omega]^{\textstyle <\!\omega}} 
\newcommand{\fseo}{\omega^{\textstyle <\!\omega}} 
\newcommand{\lh}{{\rm lh}\/} 
\newcommand{\rest}{{\mathord{\restriction}}}

\newcommand{\dom}{{\rm dom}} 
\newcommand{\rng}{{\rm rng}}
 
\newcommand{\A}{{\mathcal A}}  
  
\newcommand{\B}{{\mathcal B}}
\newcommand{\random}{{\mathbb B}}
\newcommand{\D}{{\mathbb D}} 
\newcommand{\E}{{\mathbb E}} 
\newcommand{\F}{{\mathcal F}} 
 
\newcommand{\G}{{\mathcal G}}

\newcommand{\I}{{\mathcal I}}  

\newcommand{\M}{{\mathcal M}}  
\newcommand{\m}{{\mathbb M}}

\renewcommand{\P}{{\mathcal P}}
\newcommand{\p}{{\mathbb P}}

\newcommand{\Z}{{\mathbb Z}}

\newcommand{\an}{{\bf An}}
\newcommand{\borel}{{\bf Borel}}

\newtheorem{theorem}{Theorem}[section] 
 
\newtheorem{proposition}[theorem]{Proposition} 
\newtheorem{corollary}[theorem]{Corollary} 

\theoremstyle{definition}
\newtheorem{problem}[theorem]{Problem} 
\newtheorem{definition}[theorem]{Definition}

\title{Ideals determined by some Souslin forcing notions} 
\author{Haim Judah}
\address{Dept. of Mathematics and Computer Science\\
Bar-Ilan University\\
52900 Ramat-Gan, Israel}

\author{Andrzej Ros{\l}anowski}
\address{Mathematical Institute of Wroclaw University\\
50384 Wroclaw, Poland}
\thanks{The research supported by  KBN
(Polish Committee of Scientific Research) grant 1065/P3/93/04}

\begin{document} 
\begin{abstract}
We describe a method of building ``nice'' $\sigma$--ideals from
Souslin ccc forcing notions.   
\end{abstract}
\maketitle 

\section{Introduction} 

\noindent{\bf Preliminaries:} 
By Sikorski theorem (cf \cite[\S 31]{Si64}) every ccc countably generated
complete Boolean algebra $\random$ is isomorphic to the quotient algebra
$\borel(\can)/I$ of Borel subsets of the Cantor space modulo some Borel
$\sigma$-ideal. The isomorphism can be described as follows. Let
$\mbox{ST}(\random)$ be the Stone space of the algebra $\random$, $\M$ be
the $\sigma$-ideal of meager sets of the space. Then the algebra $\random$
is isomorphic to the quotient $\mbox{BAIRE}(\mbox{ST}(\random))/\M$ of Baire
subsets of $\mbox{ST}(\random)$ modulo meager sets. Let $u_n\in\random$ be
generators of $\random$ (so $[u_n]_{\M}$ are generators of
$\mbox{BAIRE}(\mbox{ST}(\random))/\M$; elements of $\random$ are identified
with clopen subsets of $\mbox{ST}(\random)$). Define
$\phi:\mbox{ST}(\random)\longrightarrow\can$ by $\phi(x)(n)=1$ iff $x\in
u_n$ and let
$f:\borel(\can)\longrightarrow\mbox{BAIRE}(\mbox{ST}(\random))/\M: A\mapsto
[\phi^{-1}[A]]_{\M}$. Then $f$ is a $\sigma$-epimorphism of the Boolean
algebras and hence $\random\cong\borel(\can)/\mbox{Ker}(f)$. The ideal
$\mbox{Ker}(f)$ consists of all Borel sets $A\subseteq\can$ such that the
pre-image $\phi^{-1}[A]$ is meager in $\mbox{ST}(\random)$. As both the space
$\mbox{ST}(\random)$ and the ideal of meager sets of it have no nice general
description this approach has several disadvantages.  In particular it is
difficult to describe and to investigate the ideal $\mbox{Ker}(f)$. Moreover
generally it has none of properties we could expect - we should keep in mind
the ideals of meager and of null subsets of the Cantor space as ``positive''
examples here. For these reasons we will present an another approach.

One could ask why we are interested in representing a forcing notion as a
quotient Boolean algebra $\borel/\I$ for some ccc Borel $\sigma$-ideal
$\I$. There are several reasons. The first (and the main one) is that if
we assume (or prove) some additional properties of the ideal $\I$ then we
can use the well-developed machinery of $\I$-random reals (cf
\cite{Ku84}). The second reason is that we have a nice description of reals,

Borel sets etc in extensions via such algebras.

\begin{proposition}
Suppose that $\I$ is a ccc Borel $\sigma$-ideal on $\can$,
$\random=\borel/\I$ is the quotient (complete) algebra. Let $\dot{r}$ be
a $\random$-name for an element of $\can$ such that $\lbv
s\subseteq\dot{r}\rbv=[[s]]_{\I}$.
\begin{enumerate}
\item If $\tau$ is a $\random$-name for an element of $\can$ then there is
a Borel function $f:\can\longrightarrow\can$ such that
$\forces_{\random}f(\dot{r})=\tau$.
\item If $\dot{B}$ is a $\random$-name for a Borel subset of $\can$ then
there is a Borel set $A\subseteq\can\times\can$ such that
$\forces_{\B}\dot{B}=(A)_{\dot{r}}$, where $(A)_x=\{y:(x,y)\in A\}$. 
\end{enumerate}
\end{proposition}

\begin{proof} {\em 1.}\ \ \ Construct inductively Borel sets
$A_s\subseteq\can$ such that for each $s\in\fs$: $A_s=A_{s\hat{\ }0}\cup
A_{s\hat{\ }1}$, $A_{s\hat{\ }0}\cap A_{s\hat{\ }1}=\emptyset$ and $\lbv
s\subseteq\dot{r}\rbv_{\random}=[A_s]_{\I}$. Put $f(x)=y$ if
$x\in\bigcap_{n\in \omega}A_{y\rest n}$ and $f(x)=\bar{0}$ if
$x\notin\bigcap_{n\in \omega}\bigcup_{s\in 2^n}A_s$. This $f$ works.\\
{\em 2.}\ \ \ Let $\{C_n:n\in \omega\}$ enumerate the (clopen) basis of
$\can$. If $\dot{B}$ is a name for an open set then we have a name
$\dot{U}$ for a subset of $\omega$ such that
$\forces\dot{B}=\bigcup\{C_n:n\in\dot{U}\}$. Let $A_n$ be a Borel set such
that  $[A_n]_{\I}=\lbv n\in\dot{U}\rbv$ and put $A=\bigcup\{A_n\times
C_n:n\in \omega\}$. Since $[A_n]_{\I}=\lbv\dot{r}\in A_n\rbv$ we get that
$\forces (A)_{\dot{r}}=\dot{B}$. Thus we are done for open sets. Next
apply easy induction (note that $(\bigcap_{n}A_n)_x=\bigcap_{n}(A_n)_x$
and $(\neg A)_x=\neg (A)_x$).
\end{proof}

\noindent{\bf Notation:} Our notation is standard. However, in forcing
considerations we keep the convention that {\em a stronger condition is
the greater one}. $\con$ stands for the cardinality of the continuum.  A
notion of forcing $\p$ is Souslin if $\p$ is a $\Sigma^1_1$ subset
of reals and both the order $\leq_{\p}$ and the incompatibility relation
$\incomp_{\p}$ are $\Sigma^1_1$-subsets of the plane.

\section{The ideal}
For a countably generated forcing notion $\p$ we want to define an ideal
$\I_{\p}$ in Borel subsets of the reals such that the respective quotient
algebra is isomorphic to $RO(\p)$. Let $p_n\in\p$ (for $n\in \omega$) be
such that they completely generate the algebra $RO(\p)$. Let $\dot{r}$ be
a $\p$-name for a real from $\can$ such that $\lbv\dot{r}(n)=1\rbv_{RO(\p)}
=p_n$. We can define an ideal $\I^0_{\p}$ putting for a Borel set
$B\subseteq\can$  
\[B\in\I_{\p}^0\ \mbox{ if and only if }\ \forces_{\p}\dot{r}\notin B.\] 
One can easily see that this is a ccc $\sigma$-ideal containing all
singletons and such that that $\borel/\I_{\p}^0\cong RO(\p)$. If $\p$ is a
Souslin forcing notion then the ideal $\I_{\p}^0$ is absolute. 

However, there are some difficulties here. It is not so easy to calculate
the complexity of the ideal and it is not obvious how to ensure
invariance of the ideal. Moreover the embedding of the forcing notion
into the respective quotient Boolean algebra has no simple description.
For some forcing notions we will construct the ideal on the Baire space
in the way giving more possibilities to work with it.

\begin{definition}
A forcing notion $\p$ is {\em countably-1-generated} if there are
conditions $p_n\in\p$ (for $n\in \omega$) such that
\[(\forall p\in\p)(\forall q\in\p, q\incomp p)(\exists n\in \omega)
(p_n\incomp p\ \&\ p_n\compatible q).\] 
In this situation the conditions $p_n$ ($n\in \omega$) are called
$\sigma$-1-generators of the forcing notion $\p$.
\end{definition}

Clearly if $\p$ is countably-1-generated then the Boolean algebra $RO(\p)$
is countably generated and each element of $\p$ is the complement (in the
algebra $RO(\p)$) of the union of a family of generators. Thus elements of
$RO(\p)$ are unions of elements of the form $-\sum\{p_n:n\in u\}$,
$u\subseteq \omega$ but they do not have to be of this form. Many classical
ccc countably generated Boolean algebras are determined by
countably-1-generated forcing notions. The Random Algebra is determined by
the order of closed sets of positive measure in $\can$.  Clearly this order
is countably-1-generated. The Amoeba Algebra for measure, the Amoeba Algebra
for category, the Hechler forcing and the Eventually Different Real forcing
notion can be represented as countably-1-generated orders. Actually we have
no example of a ccc Souslin forcing notion (producing one real extension)
which is not of this kind.

\begin{problem}
Suppose $\p$ is a ccc Souslin forcing notion such that the algebra
$RO(\p)$ is countably generated. Can $\p$ be represented as a ccc Souslin
countably-1-generated forcing notion?
\end{problem}

In our considerations we will assume that every forcing notion is {\em
separative}, i.e. if $p,q\in\p$, $p\not\leq q$ then there is $q_0\geq q$
such that $q_0\incomp p$. This assumption can be easily avoid if we
replace (in some places) inequality in $\p$ by that in $RO(\p)$.

\begin{proposition}
\label{basis}
Suppose $\p$ is an atomless ccc countably-1-generated forcing notion.
Then there is a mapping $\pi:\fseo\longrightarrow\p$ such that 
\begin{enumerate}
\item for each $s\in\fseo$ the family $\{\pi(s\hat{\ }n):n\in \omega\}$
is a maximal antichain above $\pi(s)$,
\item $\pi(\langle\rangle)=\emptyset_{\p}$ and
\item $\rng(\pi)$ is a set of $\sigma$-1-generators for $p$.
\end{enumerate}
\end{proposition}

\begin{proof} Let $\langle p_n:n\in \p\rangle\subseteq\p$ be a sequence of  
$\sigma$-1-generators. Construct inductively infinite maximal antichains
$\A_n\subseteq\p$ such that
\begin{itemize}
\item for each $p\in\A_n$ the set $\{q\in\A_{n+1}:q\geq p\}$ is an
infinite maximal antichain above $p$, and
\item $\{q\in\A_n:p_n\leq q\}$ is a maximal antichain above $p_n$.
\end{itemize}
Use these antichains to define $\pi$ in such a way that
$\pi[\omega^{\textstyle n+1}]=\A_n$.
\end{proof}

The mapping $\pi$ given by the above proposition (i.e.  satisfying 1-3 of
\ref{basis}) will be called {\em a basis of} the forcing notion $\p$.

Note that the formula {\em a real $b$ encodes a ccc Souslin forcing notion
and (a real) $\pi$ is a basis of it} is a $\Pi^1_2$-formula; if $b$ is a
fixed code for a ccc countably-1-generated Souslin forcing notion then {\em
$\pi$ is a basis for the forcing notion coded by $b$} is $\Pi^1_1$ (see
\cite{JdSh:292}). Consequently all the notions above are suitable absolute.

Fix a ccc countable-1-generated atomless Souslin forcing notion $\p$ and a
basis $\pi:\fseo\longrightarrow\p$ for it. Let $b$ be a real encoding $\p$.

\begin{definition}
\begin{enumerate}
\item For a condition $p\in\p$ we define 
$$\phi(p)=\{x\in\baire:(\forall n\!\in\!\omega)(\pi(x\rest n)\compatible
p)\}.$$ 
\item A set $A\subseteq\baire$ is {\em $\p$-small} if there is a maximal
antichain $\A\subseteq\p$ such that
$A\cap\bigcup\{\phi(p):p\in\A\}=\emptyset$. 
\item A set $A\subseteq\baire$ is {\em $\p$-$\sigma$-small} if it can be
covered by a countable union of $\p$-small sets. The family of
$\p$-$\sigma$-small sets will be denoted by $\I_{\p}$.
\end{enumerate}
\end{definition}

\begin{proposition}
\label{smallsets}
\begin{enumerate}
\item For each $p\in\p$ the set $\phi(p)$ is closed; $\phi(\pi(s))=[s]$
for each $s\in\fseo$. If $p\leq q$ then $\phi(q)\subseteq\phi(p)$.
\item No set $\phi(p)$ (for $p\in\p$) is $\p$-$\sigma$-small, every
singleton is $\p$-small.
\item $\p$-small sets constitute an ideal, $\I_{\p}$ is a $\sigma$-ideal
of subsets of $\baire$. Every set from $\I_{\p}$ can be covered by a
$\Sigma^0_3$-set from $\I_{\p}$. 
\end{enumerate}
\end{proposition}

\begin{proof} {\em 1.}\ \ It should be clear.

\noindent{\em 2.}\ \ Let $p\in\p$ and let $\A_n\subseteq\p$ ($n\in \omega$)
be maximal antichains. We want to find $x\in\baire$ such that $(\forall
n\!\in\!\omega)(\pi(x\rest n)\compatible p)$ and $(\forall
n\!\in\!\omega)(\exists q\!\in\!\A_n)(\forall m\!\in\!\omega)(\pi(x\rest
m)\compatible q)$. Take $p_0\in\A_0$ such that $p_0\compatible p$ and find
$n_0\in \omega$ so that $\pi(\langle n_0\rangle)\compatible(p\vee p_0)$
(i.e. such that $(\exists q\!\in\!\p)(q\geq p,p_0)$).  Choose $p_1\in\A_1$
such that $p_1\compatible(p\vee p_0\vee\pi(\langle n_0\rangle))$ and let
$n_1\in \omega$ be such that $\pi(\langle n_0,n_1\rangle)\compatible(p\vee
p_1\vee p_0\vee \pi(\langle n_0\rangle))$. Continuing in this fashion we
will define $x=\langle n_0,n_1,n_2\ldots\rangle\in\baire$ which will work.

Now suppose that $x\in\baire$. To show that the singleton $\{x\}$ is
$\p$-small it is enough to prove that the set $\{p\in\p:x\notin\phi(p)\}$ is
dense in $\p$. Given $q\in\p$. Take $q_0,q_1\geq q$ such that $q_0\incomp
q_1$. There is $s\in\fseo$ with $\pi(s)\incomp q_0$ and $\pi(s)\compatible
q_1$. If $s\subseteq x$ then $x\notin\phi(q_0)$ and we are done. So suppose
that $x\rest\lh s\neq s$. Take $q_2$ stronger than both $\pi(s)$ and
$q_1$. Then $\pi(x\rest\lh s)$ and $q_2$ are incompatible and consequently
$x\notin\phi(q_2)$.
  
\noindent {\em 3.}\ \ To prove the additivity of $\p$-small sets note that
if maximal antichains $\A_i\subseteq\p$ ($i=0,1$) witness that sets
$A_i\subseteq\baire$ are $\p$-small then any maximal antichain
$\A\subseteq\p$ refining both $\A_0$ and $\A_1$ witnesses that $A_0\cup A_1$
is $\p$-small. 
\end{proof}

\begin{definition} 
Let $\dot{r}=\dot{r}_{\pi}$ be the $\p$-name for a real in $\baire$ such
that for each $s\in\fseo$ we have $\pi(s)\forces_{\p}s\subseteq
\dot{r}_{\pi}$. 

For a real $r\in\baire$ we define $G(r)=\{p\in\p:r\in\phi(p)\}$.
\end{definition}

\begin{proposition}
\label{generic}
Let $N$ be a transitive model of $\mbox{ZFC}^*$ such that $b, \pi$ and
everything relevant is in $N$. 
\begin{enumerate}
\item If $G\subseteq\p^N$ is a generic filter over $N$ then
$G(\dot{r}^G)\cap N=G$. 
\item Suppose that $x\in\baire$ is such that for any maximal antichain
$\A\subseteq\p$, $\A\in N$ we have $x\in\bigcup_{p\in\A}\phi(p)$. Then
$G(x)\cap N$ is a generic filter over $N$ and $\dot{r}^{G(x)}$=x. 
\end{enumerate} 
\end{proposition}

\begin{proof}
First note that, in $N$, b encodes a ccc Souslin forcing notion and $\pi$ is
a basis for it ($\Pi^1_2$ formulas are downward absolute for all models of
$\mbox{ZFC}^*$). Moreover if $N\models$``{\em $\A$ is a maximal antichain in
$\p$}'' then $\A$ is really a maximal antichain of $\p$. Notice that
$\P^N=\p\cap N$ and the same concerns $\incomp_{\p}$, $\leq_{\p}$.

\noindent{\em 1.}\ \ As $G$ is a filter in $\p\cap N$ we have $G\subseteq
G(\dot{r}^G)$. If $p\notin G$, $p\in\p\cap N$ then there is $s\in\fseo$ such
that $\pi(s)\incomp p$ and $\pi(s)\in G$ ($\pi$ is a basis for $\p$). 
Consequently $s\subseteq\dot{r}^G$ and $\dot{r}^G\notin\phi(p)$ (so $p\notin
G(\dot{r}^G)$. 

\noindent{\em 2.}\ \ As $x\in\phi(p)\Leftrightarrow p\in G(x)$ it is enough
to show that $G(x)\cap N$ is a filter. For this it suffices to prove that
$G(x)\cap N$ contains no pair of incompatible elements. Thus suppose that
$p_0,p_1\in\p\cap N$ are incompatible. Let $\A\in N$ be a maximal antichain
in $\p$ such that (in $N$) for each $p\in\A$
\begin{quotation}
\noindent either there is $s\in\fseo$ such that $p\geq\pi(s)$ and
$\pi(s)\incomp p_0$

\noindent or there is $s\in\fseo$ such that $p\geq\pi(s)$ and $\pi(s)\incomp
p_1$.
\end{quotation}
By the choice of $x$ we have that $x\in\phi(p)$ for some $p\in\A$. Let
$s\in\fseo$ be such that $p\geq\pi(s)$ and $\pi(s)\incomp p_0$ (or $\pi(s)
\incomp p_1$). Then $s\subseteq x$ and $x\notin\phi(p_0)$ (or $x\notin
\phi(p_1)$). Consequently either $p_0\notin G(x)$ or $p_1\notin G(x)$.
\end{proof}

\begin{proposition}
\begin{enumerate}
\item Let $B$ be a Borel subset of $\baire$. Then
\[B\notin\I_{\p}\mbox{ if and only if }(\exists p\!\in\!\p)((\phi(p)
\backslash B)\in\I_{\p}).\] 
\item The formula {\em $c$ is a code for a Borel set belonging to $\I_{\p}$}
is $\Delta^1_2$; it is absolute for all transitive models of $\mbox{ZFC}^*$.
\end{enumerate}
\end{proposition}

\begin{proof}
{\em 1.}\ \ Since $\phi(p)\notin\I_{\p}$ for any $p\in\p$ (by
\ref{smallsets}) we easily get that $(\exists p\!\in\!\p)(\phi(p)\backslash
B\in\I_{\p})$ implies $B\notin\I_{\p}$.  Suppose now that
$B\notin\I_{\p}$. Let $c$ be a real encoding the Borel set $B$. Let $N$ be a
countable transitive model of $\mbox{ZFC}^*$ such that $b,c,\pi,\ldots\in
N$. Since $B\notin I_{\p}$ we find a real $x\in B$ such that
\[x\in\bigcap\{\bigcup_{p\in\A}\phi(p):N\models\A\mbox{{\em\ is a maximal
antichain in }}\p^N\}.\] 
By \ref{generic} we get that $G=G(x)\cap N$ is a $\p^N$-generic filter over
$N$ and $\dot{r}^G=x$. As $N[G]\models\dot{r}^G\in B$ we find $p\in\p^N$
such that $N\models p\forces\dot{r}\in\sharp c$ (where $\sharp c$ stands for
the Borel set coded by $c$). We claim that
\[(\phi(p)\backslash B)\cap\bigcap\{\bigcup_{q\in\A}\phi(q):
N\models\A\mbox{{\em\ is a maximal antichain in }}\p^N\}=\emptyset.\]
Suppose not and let $y$ be a real from the intersection. As earlier we
have that $G'=G(y)\cap N$ is a $\p^N$-generic filter over $N$, $\dot{r}^{G'}
=y$. Since $p\in G(y)$ we get a contradiction to $y\notin B$.

\noindent{\em 2.}\ \ For a real $a$ let $\langle (a)_n:n\in \omega\rangle$
be the sequence of reals coded by $a$. Let $\an=\{a:\langle(a)_n:n\in
\omega\rangle\subseteq\p\mbox{ is a maximal antichain }\}$. Clearly $\an$ is
the intersection of a $\Pi^1_1$-set and a $\Sigma^1_1$-set. Now\\ 
{\em $c$ is a Borel code for a set from $\I_{\p}$} $\equiv$\\
$(\exists a)((\forall n)((a)_n\in\an)\ \&\ (\forall x\in\sharp c)(\exists
n,m)(x\in\phi(((a)_n)_m)))\ \&\ c\in\mbox{BC}$\\ 
The first part of the conjunction is $\Sigma^1_2$, the second part is
$\Pi^1_1$. Hence the formula is $\Sigma^1_2$. On the other hand, by {\em
1.},\\ 
{\em $c$ is a Borel code for a set not belonging to $\I_{\p}$} $\equiv$\\
$(\exists p\in\p)((\phi(p)\backslash\sharp c)\in\I_{\p})\ \&\
c\in\mbox{BC}$.\\ 
Easily the last formula is $\Sigma^1_2$ too. Consequently both formulas are
$\Delta^1_2$ and this fact is provable in ZFC. As $\Sigma^1_2$ formulas are
upward absolute (for models of $\mbox{ZFC}^*$) and $\Pi^1_2$ formulas are
downward absolute (for models of $\mbox{ZFC}^*$) we are done. 
\end{proof}

\begin{corollary}
$\I_{\p}$ is a Borel ccc absolute $\sigma$-ideal on $\baire$. The
quotient algebra $\borel(\baire)/\I_{\p}$ is a ccc complete Boolean
algebra. The mapping
\[\p\longrightarrow\borel(\baire)/\I_{\p}:p\mapsto[\phi(p)]_{\I_{\p}}\]
is a dense embedding (so $RO(\p)\cong\borel(\baire)/\I_{\p}$). For each
Borel code $c$: $\lbv\dot{r}\in \sharp c\rbv_{\p}=[\sharp c]_{\I_{\p}}$
(thus in particular $\I_{\p}=\I_{\p}^0$).  
\end{corollary}

\section{Invariance}
In this section we will be interested in invariant properties of ideals
$\I_{\p}$. For an ideal on $\baire$ we can consider its invariance under
permutations of $\omega$ (both as domain of sequences and as the set of
their values) as well as invariance under translations. We can equip
$\omega$ with an additive structure copied from $\Z$. Then the Baire space
becomes a group too (with the product of the addition). The invariance under
translations in this group can be captured by a more general invariance.

\begin{definition}
Let $\I$ be an ideal on $\baire$.
\begin{enumerate}
\item $\I$ is weakly index invariant if for any permutation
$P:\omega\stackrel{\rm onto}{\longrightarrow}\omega$ and a set
$A\subseteq\baire$
\[A\in\I\mbox{ if and only if }\{x\comp P:x\in A\}\in\I.\]
\item The ideal is permutation invariant if for any sequence of
permutations $P_n:\omega\stackrel{\rm onto}{\longrightarrow}\omega$
($n\in \omega$) and a set $A\subseteq\baire$
\[A\in\I\mbox{ if and only if }\{\bar{P}(x):x\in A\}\in\I,\]
where $\bar{P}(x)(n)=P_n(x(n))$.
\end{enumerate}
\end{definition}
If we want the ideal $\I_{\p}$ to be invariant we have to assume some
extra properties of the pair $(\p,\pi)$. 

\begin{definition}
\label{invariantbasis}
Let $\pi$ be a basis for a countably-1-generated Souslin forcing notion
$\p$.
\begin{enumerate}
\item The basis $\pi$ is {\em index invariant} if for every permutation
$P:\omega\stackrel{\rm onto}{\longrightarrow}\omega$ there is an
automorphism $a_P:\p\stackrel{\rm onto}{\longrightarrow}\p$ such that for
each $p\in\p$
\[\phi(a_P(p))=\{x\comp P:x\in\phi(p)\}.\]
\item The basis $\pi$ is {\em permutation invariant} if for every
sequence $\bar{P}$ of permutations $P_n:\omega\stackrel{\rm onto}
{\longrightarrow}\omega$ there is an automorphism
$a^{\bar{P}}:\p\stackrel{\rm onto}{\longrightarrow}\p$ such that for every
$p\in\p$ 
\[\phi(a^{\bar{P}}(p))=\{\bar{P}(x):x\in\phi(p)\}.\]
\end{enumerate}
\end{definition}

Directly from the definition we can conclude the following observation.
\begin{proposition}
\label{indexpermutation}
If the basis $\pi$ for $\p$ is index invariant then the ideal
$\I_{\p}$ is weakly index invariant. If the basis is permutation
invariant then the ideal is permutation invariant. 
\end{proposition}

\begin{definition}
\begin{enumerate}
\item For $K\in\iso$ let $\mu_K:\omega\stackrel{\rm
onto}{\longrightarrow}K$ be the increasing enumeration.
\item An ideal $\I$ on $\baire$ is {\em injective} if for each set
$K\in\iso$ and a set $A\subseteq \omega^{\textstyle K}$
\[\{x\in\baire:x\rest K\in A\}\in\I\ \mbox{ if and only if }\
\{x\comp\mu_K: x\in A\}\in\I.\]
\item The ideal is {\em index invariant} if it is weakly invariant and
injective. 
\end{enumerate}
\end{definition}

\begin{definition}
A basis $\pi$ for a forcing notion $\p$ is {\em productive}
if for every $K\in\iso$ there is a complete embedding
$i_K:\p\longrightarrow\p$ such that 
\[(\forall p\in\p)(\phi(i_K(p))=\phi_K(p)\times \omega^{\textstyle
\omega\backslash K}),\]
where $\phi_K(p)=\{x\comp\mu_K^{-1}: x\in\phi(p)\}$.
\end{definition}

\begin{proposition}
Suppose that $\pi$x is a productive basis for $\p$. Assume that each
$\Pi^1_1$ subset of $\baire$ is either $\p$-$\sigma$-small or
$\I_{\p}$-almost contains a set $\phi(p)$ for some $p\in\p$. Then the
ideal $\I_{\p}$ is injective
\end{proposition}

\begin{proof}
Let $K\in\iso$ and $A\in \omega^{\textstyle K}$.

Suppose that $\{x\comp\mu_K: x\in A\}$ is $\p$-small. Let $\A\subseteq\p$ be
a maximal antichain in $\p$ such that $(\forall x\in A)(\forall
p\in\A)(x\comp\mu_K\notin\phi(p))$. Look at $i_K[\A]$ (where
$i_K:\p\longrightarrow\p$ is the embedding given by productivity of
$\pi$). It is a maximal antichain in $\p$. Suppose now that $x\in\baire$,
$x\rest K\in A$ and $p\in\A$. Then $x\comp\mu_K=(x\rest K)\comp\mu_K\notin
\phi(p)$ and hence $x\rest K=(x\comp\mu_K)\comp\m_K^{-1}\notin\phi_K(p)$. As
$\phi(i_K(p))=\phi_K(p)\times \omega^{\textstyle \omega\backslash K}$ we
conclude $x\notin\phi(i_K(p))$. Consequently $i_K[\A]$ witnesses that
$\{x\in\baire: x\rest K\in A\}$ is $\p$-small. \\ Now we can easily conclude
that $\{x\comp\mu_K:x\in A\}\in\I_{\p}$ implies $\{x\in\baire:x\rest K\in
A\}\in\I_{\p}$.

Assume now that $\{x\in\baire: x\rest K\in A\}\in\I_{\p}$. Let $B\subseteq
\baire$ such that $B\in\I_{\p}$ and $\{x\in\baire:x\rest K\in A\}\subseteq
B$. Let $B^*=\{c\in \omega^{\textstyle K}: (\forall x\in\baire)(c\subseteq
x\Rightarrow x\in B)\}$. Clearly $B^*$ is a $\Pi^1_1$ subset of
$\omega^{\textstyle K}$ and $A\subseteq B^*$. Thus if we prove that $\{c
\comp\mu_K:c\in B^*\}\in\I_{\p}$ then we will have $\{c\comp\mu_K:c\in A\}
\in\I_{\p}$ and the proposition will be proved.\\ 
Suppose $\{c\comp\mu_K:c\in B^*\}\notin\I_{\p}$. Since it is a
$\Pi^1_1$-subset of $\baire$ we find $p\in\p$ such that $\phi(p)\backslash
\{c\comp\mu_K:c\in B^*\}\in\I_{\p}$. By the first part we get
\[\{x\in\baire: x\comp\mu_K\in\phi(p)\backslash\{c\comp\mu_K:c\in
B^*\}\}\in\I_{\p}.\]
Let $x\in\phi(i_K(p))=\phi_K(p)\times \omega^{\textstyle \omega\backslash
K}$ be such that $x\comp\mu_K\notin\phi(p)\backslash\{c\comp\mu_K:c\in
B^*\}$. As $x\rest K\in\phi_K(p)$ we have $x\comp\mu_K\in\phi(p)$ and
hence $x\comp\mu_K=c\comp\mu_K$ for some $c\in B^*$. This means
$c\subseteq x$ and by the definition of $B^*$ this implies $x\in B$. Thus
we have proved that $\phi(i_K(p))\backslash B\in\I_{\p}$ what implies
$B\notin\I_{\p}$ - a contradiction. 
\end{proof}

\section{Baire Property}
\begin{definition}
A family $\F$ of subsets of $\baire$ is {\em a category base on $\baire$}
if $|\F|=\con$, $\bigcup\F=\baire$ and\\ 
for each subfamily $\G\subseteq\F$ of disjoint sets, $|\G|<\con$ and each
$A\in\F$
\begin{quotation}
\noindent if $(\exists B\!\in\!\F)(B\subseteq A\cap\bigcup\G)$ then
$(\exists B\!\in\!\F)(\exists C\!\in\!\G)(B\subseteq C\cap A)$, and

\noindent if there is no $B\in\F$ contained in $A\cap\bigcup\G$ then
$(\exists B\!\in\!\F)(B\subseteq A\backslash\bigcup\G)$.
\end{quotation}
\end{definition}

\begin{definition}
Let $\F$ be a category base on $\baire$ and let $A\subseteq\baire$.
\begin{enumerate}
\item $A$ is $\F$-singular if $(\forall B\!\in\!\F)(\exists
c\!\in\!\F)(C\subseteq B\backslash A)$.
\item $A$ is $\F$-meager if it can be covered by a countable union of
$\F$-singular sets. 
\item $A$ has $\F$-Baire property if for every $B\in\F$ there is $C\in\F$
such that either $C\cap A$ is $\F$-meager or $C\backslash A$ is $\F$-meager.
\end{enumerate}
\end{definition}

\begin{theorem}[Marczewski, Morgan]
Let $\F$ be a category base on $\baire$. Then $\F$-meager sets constitute
a $\sigma$-ideal on $\baire$. Sets with the $\F$-Baire property form a
$\sigma$-field which is closed under the Souslin operation $\A$. 
\end{theorem}

Suppose that $\pi$ is a basis for a ccc Souslin forcing $\p$. Assume that
\[(\forall p,q\in\p)(p\incomp q\ \mbox{ if and only if }\
\phi(p)\cap\phi(q)=\emptyset).\] 

\begin{proposition}
If $\pi,\p$ are as above then the family $\F_{\p}=\{\phi(p): p\in\p\}$ is
a category base. The family of $\F_{\p}$-singular sets is the family of
$\p$-small subsets of $\baire$, $\F_{\p}$-meager sets agree with
$\P$-$\sigma$-small sets. 
\end{proposition}

\begin{corollary}
\label{BaireInjective}
Let $\pi,\p$ be as above.
\begin{enumerate}
\item Sets with $\F_{\p}$-Baire property constitute a $\sigma$-field of
subsets of $\baire$. This $\sigma$-field is closed under the Souslin
operation $\A$ and contains all Borel sets (and hence it contains both
$\Sigma^1_1$ and $\Pi^1_1$ sets).
\item If $A\in\baire$ is a $\Pi^1_1$ set then either $A\in\I_{\p}$ or
$(\exists p\in\p)(\phi(p)\backslash A\in\I_{\p})$.
\item If additionally $\pi$ is productive then the ideal $\I_{\p}$ is
injective. 
\end{enumerate}
\end{corollary}

\section{Two examples}
The Eventually Different Real forcing notion $\E$ can be represented as a
family of pairs $(z,F)$ where $z$ is a finite function $z:\dom(z)
\longrightarrow\omega$, $\dom(z)\subseteq \omega$ and $F:\omega
\longrightarrow\fsuo$ is such that $\max\{|F(n):n\in\omega\}<\omega$. The
pairs are ordered by   
\begin{quotation}
\noindent $(z,F)\leq (z',F')$ if and only if

\noindent $z\subseteq z'$, $F\subseteq F'$ and $(\forall k\!\in\!(\dom
z\backslash\dom z'))(z'(k)\notin F(k))$.
\end{quotation}
It is well known (see \cite{Mi81}) that this order is a ccc Souslin
forcing notion (even with the small modification we have introduced). For
$s\in\fseo$ we put $\pi_{\E}(s)=(s,F_{\emptyset})$, where
$F_{\emptyset}(n)=\emptyset$.  We claim that $\pi_{\E}$ is a basis for
$\E$. Conditions 1, 2 of \ref{basis} are clearly satisfied. For the third
condition note that
\begin{quotation}
\noindent $(z,F), (z',F')\in\E$ are incompatible if and only if

\noindent either $z,z'$ are incompatible (as functions)

\noindent or $(\exists n\!\in\!(\dom z'\backslash\dom z))(z'(n)\in F(n))$

\noindent or $(\exists n\!\in\!(\dom z\backslash\dom z'))(z(n)\in F'(n))$
\end{quotation}
Suppose now that $(z,F)\incomp (z',F')$. We may assume that $\dom z'=m$ for
some $m\in \omega$. If $z,z'$ are incompatible then $\pi_{\E}(z')\incomp
(z,F)$ and clearly $\pi_{\E}(z')\compatible (z',F')$. If for some $n\in
m\backslash\dom z$ we have $z'(n)\in F(n)$ then $\pi_{\E}(z')$ is
incompatible with $(z,F)$ and compatible with $(z',F')$.  So suppose that
$z(n)\in F'(n)$ for some $n\in\dom z\backslash m$. Let $z''\in
\omega^{\textstyle n+1}$ be such that $z'\subseteq z''$ and $z''(k)>\max
F(k)$ for each $k\in [m,n]$. Then $z'',z$ are incompatible (so
$\pi_{\E}(z'')\incomp (z,F)$) but $\pi_{\E}(z'')\compatible (z',F')$.

Note that $\phi(z,F)=\{x\in\baire: z\subseteq F\ \&\ (\forall n\notin\dom
z)(x(n)\notin F(n))\}$. 

\begin{proposition}
\label{EvDiffReal}
The basis $\pi_{\E}$ is index invariant permutation invariant basis
productive basis. Moreover for each $p,q\in\E$
\[\phi(p)\cap\phi(q)=\emptyset\ \mbox{ if and only if }\ p\incomp q.\]
\end{proposition}

\begin{proof}
We want to show the invariance of the basis $\pi_{\E}$. For permutations
$P,P_n:\omega\stackrel{\rm onto}{\longrightarrow}\omega$ we define 
\[a_P:\E\longrightarrow\E:(z,F)\mapsto(z\comp P, F\comp P),\] 
\[a^{\bar{P}}:\E\longrightarrow\E:(z,F)\mapsto(\bar{P}(z), \bar{P}(F))\] 
(where $z\comp P$ is defined on $P^{-1}[\dom z]$, $\bar{P}(z)$ is
defined on $\dom z$, $\bar{P}(z)(n)=P_n(z(n))$ and
$\bar{P}(F)(n)=P_n[F(n)]$). We claim that  
$a_P$, $a^{\bar{P}}$ have the properties required in definition
\ref{invariantbasis}. Clearly both are automorphisms of $\E$.
Moreover,\\
$x\in\phi(a_P(z,F))\ \equiv\ x\in\phi(z\comp P, F\comp P)\ \equiv$\\
$(\forall k\!\in\! P^{-1}[\dom z])(x(k)=z(P(k)))\ \&\ (\forall
m\!\notin\! P^{-1}[\dom z])(x(m)\notin F(P(m)))\ \equiv$\\ 
$(\forall k\!\in\!\dom z)(x(P^{-1}(k))=z(k))\ \&\ (\forall
m\!\notin\!\dom z)(x(P^{-1}(m))\notin F(m))\ \equiv$\\
$x\comp P^{-1}\in\phi(z,F)$.\\
Hence $\phi(a_P(z,F))=\{x\comp P: x\in\phi(z,F)\}$. Similarly\\
$x\in\phi(a^{\bar{P}}(z,F))\ \equiv\ x\in \phi(\bar{P}(z),\bar{P}(F))\
\equiv$\\ 
$(\forall k\!\in\!\dom z)(x(k)=P_k(z(k)))\ \&\ (\forall
m\!\notin\!\dom z)(x(m)\notin P_m[F(m)])\ \equiv$\\
$(\forall k\!\in\!\dom z)((P^{-1}_k\comp x)(k)=z(k))\ \&\ (\forall
m\!\notin\!\dom z)((P^{-1}_m\comp x)(m)\notin F(m))\ \equiv$\\ 
$\bar{P}^{-1}\comp x\in\phi(z,F).$\\
Thus $\phi(a^{\bar{P}}(z,F))=\{\bar{P}\comp x:x\in\phi(z,F)\}$. 

The basis $\pi_{\E}$ is productive. If for $K\in\iso$ and $(z,F)\in\E$ we
define $i_K(z,F)=(z\comp\mu_K^{-1}, (F\comp\mu_K^{-1})^*)$ (where $^*$c
means that we extend the function putting $\emptyset$ whenever it is not
defined) then $i_K$ is a complete embedding of $\E$ into $\E$ with the
required property.
\end{proof}

Applying \ref{indexpermutation}, \ref{BaireInjective} and
\ref{EvDiffReal} we get the following result.
\begin{theorem}
\label{EDR}
The ideal $\I_{\E}$ determined by the Eventually Different Real forcing
notion is a Borel ccc index invariant, permutation invariant absolute
$\sigma$-ideal on the Baire space $\baire$.  The quotient algebra
$\borel/\I_{\E}$ is neither the Cohen algebra nor the Solovay algebra.
\end{theorem}

The Hechler forcing notion $\D$ we represent as the set of pairs $(z,F)$
such that $F:\omega\longrightarrow\fsuo$ and
$z:\dom(z)\longrightarrow\omega$, $\dom(z)\in\fsuo$. The order is given by
\begin{quotation}
\noindent $(z,F)\leq (z',F')$ if and only if

\noindent $z\subseteq z'$, $F\subseteq F'$ and $(\forall k\!\in\!(\dom
z\backslash\dom z'))(z'(k)\notin F(k))$.
\end{quotation}

It is a ccc Souslin forcing notion. If we put $\pi_{\D}(s)=(s,
F_{\emptyset})$, for $s\in\fseo$ then, similarly as for $\pi_{\E}$, one can
show that $\pi_{\D}$ is a basis for $\D$. Moreover the same arguments as in
the case of the Eventually Different Real forcing show the following.

\begin{theorem}
\label{Hechler}
The ideal $\I_{\D}$ is a index invariant permutation invariant ccc Borel 
$\sigma$-ideal on the Baire space. This ideal is absolute. The forcing
with the quotient algebra $\borel/\I_{\p}$ adds a dominating real
\end{theorem}

Kunen (cf \cite{Ku84}) asked if there exists ideals with properties
mentioned in \ref{EDR}, \ref{Hechler} on the Cantor space $\can$. The ideals
$\I_{\E}$, $I_{\D}$ are not the right examples but these results show that a
large part of $\I$-random reals machinery can be allied to them as well.
\medskip

\noindent{\bf Added in May 1999:}\qquad These notes, in a slightly revised
form, were incorporated to Bartoszy\'nski and Judah \cite[\S 3.7]{BaJu95}. 

The problem of Kunen mentioned above was fully solved in Ros{\l}anowski and
Shelah \cite{RoSh:628}. The investigations of ccc ideals determined by
Souslin ccc forcing notions in a way presented here are continued in
Ros{\l}anowski and Shelah \cite[\S 3]{RoSh:672}.


\end{document}